\documentstyle[12pt]{article}

\newtheorem{Thm}{Theorem}[section]

\newtheorem{Lem}[Thm]{Lemma}

\newcommand{\subk}{{_{\mbox {\scriptsize $k$}}}}

\def \N{\rm {\bf N}}
\def \R{\rm {\bf R}}

\begin{document}

\title{A counterexample to a question of\\ R. Haydon, E. Odell and
H. Rosenthal}

\author{G. Androulakis \thanks{This work is part of the author's
Ph.D. thesis which was completed at the University of Texas at Austin
under the supervision of Professor H. Rosenthal.}}

\date{October 13, 1996}
\maketitle

\noindent
{\bf Abstract:}
We give an example of a compact metric space $K$, an open
dense subset $U$ of $K$, and a sequence
$(f_n)$  in $C(K)$ which is pointwise convergent to a non-continuous
function on $K$, such that for
every $u \in U$ there exists $n \in \N$ with $f_n(u)=f_m(u)$ for all $m
\geq n$, yet $(f_n)$ is equivalent to the unit vector basis of the
James  quasi-reflexive space of
order 1. Thus $c_0$ does not embed isomorphically in the closed linear
span $[f_n]$ of $(f_n)$. This answers in negative a question asked by H. Haydon,
E. Odell and H. Rosenthal.

\bigskip
\noindent
\section{Introduction} \label{S:intro}
A result of J. Elton \cite{E}, which was also proved later by
R. Haydon, E. Odell and H. Rosenthal \cite{HOR}, states that if $K$ is
a compact metric space, and $(f_n)$ is a uniformly bounded sequence in
$C(K)$ such that 
\[ \sum_{n=1}^\infty | f_{n+1}(k)-f_n(k)| < \infty, \; \forall k \in K
\]
and the pointwise limit of $(f_n)$ on $K$ is a non-continuous
function, then $c_0$ embeds isomorphically in the closed linear span
$[f_n]$ of $(f_n)$.
Thus the following question was
naturally raised by R. Haydon, E. Odell and H. Rosenthal:

\bigskip
\noindent
{\bf Question 4.7 in \cite{HOR}:} Let $K$ be a compact metric space, $R$ be a residual
subset of $K$ (i.e. $K \backslash R$ is a first category set), and $(f_n)$ be a
sequence in $C(K)$ which converges pointwise on $K$ to a
non-continuous function, and
\[ \sum_{n=1}^\infty | f_{n+1}(r)-f_n(r)|< \infty, \mbox{ for all }r
\in R. \]
Does $c_0$ embed in the closed linear span $[f_n]$ of $(f_n)$?

\bigskip
\noindent
We will construct $K$ a compact metric space, $U$ an open dense subset
of $K$ and a  sequence $(g_n) \subset C(K)$ such that 
\begin{itemize}
\item[(a)] $(\sum_{i=1}^n
g_i)_n$ is a uniformly bounded and pointwise convergent sequence on $K$ to
a  non-continuous function;
\item[(b)] For every $u
\in U$ there exists $n \in \N$ such that $g_m(u)=0$ for every $m
\geq n$;
\item[(c)] $[g_n]$ is isomorphic to the James quasi-reflexive of order
  1 space $J$.
\end{itemize}
Since, of course, $c_0$ does not embed isomorphically in $J$, this answers in the
negative Question 4.7 of \cite{HOR}. Our construction is very
elementary and explicit even though a shorter proof of the existence
of a counterexample to Question 4.7 of \cite{HOR} can be given along
similar lines using more advanced machinery.

\section{The construction} \label{S:construction}
\noindent
We recall the definition of the James space $J$ and some simple
facts. Let $c_{00}$ denote the finitely supported sequences of real
numbers.  For  $(x_n) \in c_{00}$ we define 
\begin{eqnarray*}
& & \ \| (x_n) \|_J = \sup \{ [ x_{p_1}^2 + (x_{p_2}-x_{p_1})^2+ \cdots +
(x_{p_k}-x_{p_{k-1}})^2 ]^{1/2}: \hfill \\
& & \hskip 2in k \in \N, 1 \leq p_1 <p_2< \cdots <
p_{k-1}<p_k \}. \ 
\end{eqnarray*}
Then the James space $J$ is the completion of $(c_{00}, \| \cdot
\|_J)$. If $(e_n)$ is the unit vector basis of $c_{00}$, then $(e_n)$
becomes the unit vector basis of $J$, which is monotone and
shrinking. Also, $(\sum_{i=1}^n e_n)_n$ is a weak-Cauchy
sequence which is not weakly convergent in $J$. If $(a_n) \in c_0$
such that $(a_n)$  is a monotone
sequence of real numbers (i.e. non-increasing, or non-decreasing) then
$\| (a_n) \|_J = |a_1|$ (this is because if $a,b \in \R$ with $ab \geq
0$, then $a^2 + b^2 \leq (a+b)^2$). 

\bigskip
\noindent
{\bf Notation:} For $(a_n), (b_n) \in c_{00}$, we define $(a_n) \star
(b_n) \in c_{00}$, by 
\[ (a_n) \star (b_n) = (a_n b_n). \]

\begin{Lem} \label{L:inequalityj}
For $(a_n), (b_n) \in c_{00}$ we have
\[ \| (a_n) \star (b_n) \|_J \leq \| (a_n) \|_J \| (b_n) \|_\infty +
\| (a_n) \|_\infty \| (b_n) \|_J. \]
\end{Lem}

\noindent
{\bf Proof} For some $k \in \N$ and some finite sequence of positive
integers  $1 \leq p_1 < p_2< \cdots p_k$  we have:
\begin{eqnarray*}
\| (a_n) \star (b_n) \|_J &=& [ (a_{p_1}b_{p_1})^2 +
(a_{p_2}b_{p_2}-a_{p_1}b_{p_1})^2+ \cdots +
(a_{p_k}b_{p_k}-a_{p_{k-1}}b_{p_{k-1}})^2]^{1/2} \\
&=& [(a_{p_1}b_{p_1})^2 +
(a_{p_2}(b_{p_2}-b_{p_1})+(a_{p_2}-a_{p_1})b_{p_1})^2+ \cdots +\\
& & (a_{p_k}(b_{p_k}-b_{p_{k-1}})+(a_{p_k}-a_{p_{k-1}})b_{p_{k-1}})^2]^{1/2}.
\end{eqnarray*}
Therefore by the triangle inequality in $\ell_2$ we have that
\begin{eqnarray*}
\| (a_n) \star (b_n) \|_J & \leq & [a_{p_1}^2 b_{p_1}^2 +
(a_{p_2}-a_{p_1})^2 b_{p_1}^2 + \cdots + (a_{p_k}-a_{p_{k-1}})^2
b_{p_{k-1}}^2]^{1/2} + \\
& & [a_{p_2}^2(b_{p_2}-b_{p_1})^2 + \cdots +
a_{p_k}^2 (b_{p_k}-b_{p_{k-1}})^2]^{1/2} \\
& \leq & [a_{p_1}^2 + (a_{p_2}-a_{p_1})^2 + \cdots +
(a_{p_k}-a_{p_{k-1}})^2]^{1/2} \| (b_n) \|_\infty + \\
& & \| (a_n) \|_\infty
[ (b_{p_2}-b_{p_1})^2 + \cdots + (b_{p_k}-b_{p_{k-1}})^2]^{1/2} \\
& \leq & \| (a_n) \|_J \| (b_n) \|_\infty + \| (a_n) \|_\infty \|
(b_n) \|_J
\end{eqnarray*}
which finishes the proof of the lemma. \hfill $\Box$

\bigskip
\noindent
Now we are ready to see the counterexample. Let $K := \{ (a,b) \in
\R^2: 0 \leq a \leq 1, \/ 0 \leq b \leq 1 \}$. Since $C[0,1]$ is
universal for the class of separable spaces, there exists a sequence $(f_n)
\subset C[0,1]$, and $M>0$ such that $(f_n)$ is $M$-equivalent to the
unit vector basis of $J$. For $n \in \N$ set $K_n := \{ (a,b) \in
\R^2: 0 \leq a \leq 1, \/ 1/2^n \leq b \leq 1 \}$, $ R_n := \{ (a,b)
\in \R^2: 0 \leq a \leq 1, \/ 1/2^n < b \leq 1 \}$, $L_n := \{
(a,b) \in \R^2 : 0 \leq a \leq 1, \/ b= 1/2^n \}$ and $L:= \{
(a,0): 0 \leq a \leq 1 \}$. Now, for $n \in \N$
define $g_n:K \rightarrow \R$ by 
\begin{itemize}
\item $g_n \mid K_n \equiv 0$,
\item for every $0 \leq a \leq 1$, $g_n$ restricted on the segment
  connecting the points $(a,1/2^n)$ and $(a,0)$, is linear, 
\item $g_n \mid L \equiv f_n$.
\item $g_n$ is continuous,
\end{itemize}
We will show that $(g_n)$ is equivalent to the unit vector basis
$(e_i)$ of the
James space. This will imply that $(\sum_{i=1}^n g_i)_n$ is a 
weak Cauchy sequence which is not weakly convergent, which will finish
the proof.  Let $n \in \N$ and
$(\lambda_i)_{i=1}^n \subset \R$. We want to estimate $\| \lambda_1
g_1 + \cdots + \lambda_n g_n \|_\infty$. For $(a,b), (c,d) \in K$, let
$[(a,b),(c,d)]$ denote the linear segment connecting the points
$(a,b)$ and $(c,d)$. For every $0 \leq a \leq 1$ we have that 
\begin{itemize}
\item $(\lambda_1 g_1 + \cdots +\lambda_n g_n) \mid [(a,1),(a,1/2)]
  \equiv 0$,
\item $(\lambda_1 g_1 + \cdots +\lambda_n g_n) \mid [(a,1/2^i),
  (a,1/2^{i+1})]$ is linear, for every $i=1, \ldots , n-1$,
\item $(\lambda_1 g_1 + \cdots + \lambda_n g_n) \mid [(a,1/2^n), (a,0)]$
  is linear,
\item $\lambda_1 g_1 + \cdots + \lambda_n g_n$ is continuous on $K$. 
\end{itemize}
Therefore we obtain:
\begin{eqnarray*}
& & \| \lambda_1 g_1 + \cdots + \lambda_n g_n \|_\infty \\
& & \hskip .5in = \max_{2 \leq
  k \leq n} \| (\lambda_1 g_1 + \cdots + \lambda_n g_n) \mid L_\subk
\|_\infty \vee \| (\lambda_1 g_1 + \cdots + \lambda_n g_n) \mid L
\|_\infty \\
& & \hskip .5in = \max_{2 \leq k \leq n} \| (\lambda_1 g_1 + \cdots
\lambda_{k-1}  g_{k-1})
\mid L_\subk \|_\infty \vee \| \lambda_1 f_1 + \cdots + \lambda_n f_n
\|_\infty.
\end{eqnarray*}
Therefore we obtain immediately the lower estimate:
\begin{eqnarray*}
\| \lambda_1 g_1 + \cdots + \lambda_n  g_n \|_\infty & \geq & \| \lambda_1
f_1 + \cdots + \lambda_n f_n \|_\infty \\
& \geq & \frac{1}{M} \| \lambda_1 e_1 + \cdots + \lambda_n e_n \|_J.
\end{eqnarray*}
For the upper estimate we need to estimate $\| (\lambda_1 g_1 + \cdots +
\lambda_n g_n \mid L_\subk) \|_\infty$ for $2 \leq k \leq n$. Note that for $0
\leq a \leq 1$ and $2 \leq k \leq n$ we have that 
\begin{eqnarray*}
& & ( \lambda_1 g_1 + \cdots + \lambda_n g_n)(a, 1/2^k)\\
& & \hskip .5in =\lambda_1 \frac{ \frac{1}{2}-
  \frac{1}{2^k}}{\frac{1}{2}}f_1 (a) +
\lambda_2 \frac{\frac{1}{2^2}- \frac{1}{2^k}}{\frac{1}{2^2}}f_2(a) +
\cdots + \lambda_{k-1} \frac{ \frac{1}{2^{k-1}} - \frac{1}{2^k}}{
  \frac{1}{2^{k-1}}} f_{k-1}(a)\\
& & \hskip .5in = \lambda_1 \frac{2^{k-1}-1}{2^{k-1}}f_1(a)+ \lambda_2
\frac{2^{k-2}-1}{2^{k-2}} f_2(a)+ \cdots + \lambda_{k-1}
\frac{2-1}{2}f_{k-1}(a). 
\end{eqnarray*}
Therefore we have that 
\begin{eqnarray*}
& & \| \lambda_1 g_1 + \cdots + \lambda_{k-1}g_{k-1} \mid L_\subk
\|_\infty \\
& & \hskip .5in =\| \lambda_1 \frac{2^{k-1}-1}{2^{k-1}} f_1 + \lambda_2
\frac{2^{k-2}-1}{2^{k-2}} f_2 + \cdots + \lambda_{k-1} \frac{2-1}{2}
f_{k-1}  \|_\infty \\
& &  \hskip .5in \leq M \|\lambda_1 \frac{2^{k-1}-1}{2^{k-1}} e_1 + \lambda_2
\frac{2^{k-2}-1}{2^{k-2}} e_2 + \cdots + \lambda_{k-1} \frac{2-1}{2}
e_{k-1}  \|_J\\ 
& &  \hskip .5in  =M \| (\lambda_1, \lambda_2, \cdots, \lambda_{k-1},
0, \ldots)\\ 
& & \hskip 1in \star (\frac{2^{k-1}-1}{2^{k-1}},
\frac{2^{k-2}-1}{2^{k-2}},  \ldots ,
\frac{2-1}{2}, 0, \ldots) \|_J \\
& & \hskip .5in \leq M \| \lambda_1 e_1 + \cdots +
\lambda_{k-1}e_{k-1} \|_J  \cdot 1 \\
& & \hskip .5in + M \|
(\lambda_i)_{i=1}^{k-1} \|_\infty \| ( \frac{2^{k-1}-1}{2^{k-1}}, \ldots,
\frac{2-1}{2},0, \ldots) \|_J (\mbox{by Lemma
    \ref{L:inequalityj}}) \\
& & \hskip .5in \leq M \| \lambda_1 e_1 + \cdots +
\lambda_{k-1}e_{k-1} \|_J +  M \|
(\lambda_i) \|_\infty \frac{2^{k-1}-1}{2^{k-1}} \/ (\mbox{since the}\\
& & \hskip .5in \mbox{sequence }
  (\frac{2^{k-1}-1}{2^{k-1}}, \frac{2^{k-2}-1}{2^{k-2}}, \ldots ,
  \frac{2-1}{2}, 0, \ldots, ) \mbox{ is decreasing})\\
& & \hskip .5in \leq 2 M \| \lambda_1 e_1 + \cdots + \lambda_{k-1}
e_{k-1} \|_J \/ (\mbox{since } \|(\lambda_i)_{i=1}^{k-1} \|_\infty \leq \|
  (\lambda_i)_{i=1}^{k-1}   \|_J). 
\end{eqnarray*}
Also, since $\| \lambda_1 f_1 + \cdots + \lambda_n f_n \|_J \leq M \|
\lambda_1 e_1 + \cdots + \lambda_n e_n \|_J$, we obtain that
\[ \| \lambda_1 g_1 + \cdots + \lambda_n g_n \|_\infty \leq 2 M \|
\lambda_1 e_1 + \cdots + \lambda_n e_n \|_J. \]
This finishes the proof. \hfill $\Box$

\noindent
G. Androulakis, Math. Sci. Bldg., University of Missouri-Columbia, 
Columbia, MO 65211-0001\\
e-mail: giorgis@math.missouri.edu 

\end{document}